\documentclass[12pt]{amsart}
\usepackage{amssymb}
\usepackage{latexsym}
\usepackage{amsmath}

\newtheorem{Theorem}{Theorem}[section]
\newtheorem{Lemma}[Theorem]{Lemma}
\newtheorem{Corollary}[Theorem]{Corollary}
\newtheorem{Proposition}[Theorem]{Proposition}

\newenvironment{pf}{\medskip\par\noindent{\bf Proof\/}.}{\hfill
$\Box$\medskip}

\begin{document}










\title{Symplectic Graphs and Their Automorphisms}
\author{Zhongming Tang and Zhe-xian Wan}

\begin{abstract}

A new family of strongly regular graphs, called the general
symplectic graphs $Sp(2\nu, q)$, associated with nonsingular
alternate matrices is introduced. Their parameters as strongly
regular graphs, their chromatic numbers as well as their groups of
graph automorphisms are determined.

\end{abstract}

\thanks{Both Authors are Supported by the National Natural Science Foundation of China}
\date{}
\maketitle

\section{Introduction}
Let $\mathbb{F}_q$ be a finite field of any characteristic and
$\nu\geq 1$ an integer. Let
$$
\mathbb{F}_q^{(2\nu)}=\{(a_1,\ldots,a_{2\nu}):a_i\in
\mathbb{F}_q,i=1,\ldots,2\nu\}.
$$
be the $2\nu$-dimensional row vector space over $\mathbb{F}_q$.
For any $\alpha_1,\ldots,\alpha_n\in \mathbb{F}_q^{(2\nu)}$, we
denote the subspace of $\mathbb{F}_q^{(2\nu)}$ generated by
$\alpha_1,\ldots,\alpha_n$ by $[\alpha_1,\ldots,\alpha_n]$. Thus,
if $\alpha=(a_1,\ldots,a_{2\nu})\not=0\in \mathbb{F}_q^{(2\nu)}$
then $[\alpha]$, which is also denoted by $[a_1,\ldots,a_{2\nu}]$,
is a one dimensional subspace of $\mathbb{F}_q^{(2\nu)}$ and
$[\alpha]=[k\alpha]$ for any $k\in
\mathbb{F}_q^*=\mathbb{F}_q\setminus\{0\}$.

Let $K$ be a $2\nu\times 2\nu$ nonsingular alternate matrix over
$\mathbb{F}_q$. The {\em symplectic graph} relative to $K$ over
$\mathbb{F}_q$ is the graph with the set of one dimensional
subspaces of $\mathbb{F}_q^{(2\nu)}$ as its vertex set and with
the adjacency defined by
$$
[\alpha]\sim[\beta]\mbox{ if and only if }\alpha K^t\!\beta\not=0,
\mbox{ for any }\alpha\not=0, \beta\not=0 \in
\mathbb{F}_q^{(2\nu)},
$$
where $[\alpha]\sim[\beta]$ means that $[\alpha]$ and $[\beta]$
are adjacent. Since any two $2\nu\times 2\nu$ nonsingular
alternate matrices over $\mathbb{F}_q$ are cogredient, any two
symplectic graphs relative to two different $2\nu\times 2\nu$
nonsingular alternate matrices over $\mathbb{F}_q$ are isomorphic.
Thus we can assume that
$$
K=\left(
\begin{array}{cccc}
\begin{array}{cc}
0&1\\-1&0
\end{array}
&&&\\
&\begin{array}{cc} 0&1\\-1&0
\end{array}
&&\\
&&\ddots&\\
&&&\begin{array}{cc} 0&1\\-1&0
\end{array}
\end{array}
\right)_{2\nu\times 2\nu}.
$$
and consider only the symplectic graph relative to the above $K$
over $\mathbb{F}_q$, which will be denoted by  $Sp(2\nu,q)$.

When $q=2$, the special case $Sp(2\nu,2)$ of the graph
$Sp(2\nu,q)$ was studied previously by Rotman \cite{R}, Rotman and
Weichsel \cite{RW}, Godsil and Royle \cite{GR1,GR}, etc. In the
present paper we study the general case $Sp(2\nu,q)$. In Section
2, we show that $Sp(2\nu,q)$ is strongly regular and compute its
parameters. We also prove that the chromatic number of
$Sp(2\nu,q)$ is $q^{\nu}+1$. Section 3 is devoted to discuss the
group of automorphisms $\mbox{Aut}(Sp(2\nu,q))$ of the graph. The
structure of this group depends on $q$ and $\nu$. When $q=2$,
$\mbox{Aut}(Sp(2\nu,2))$ is isomorphic to the symplectic group of
degree $2\nu$ over $\mathbb{F}_2$. When $q>2$,
$\mbox{Aut}(Sp(2\nu,q))$ is the product of two subgroups which are
identified clearly (cf. Theorem 3.4).

\medskip
\section{Strongly Regularity and Chromatic Numbers of Symplectic Graphs}

For any subspace $V$ of $\mathbb{F}_q^{(2\nu)}$, we denote the
subspace of $\mathbb{F}_q^{(2\nu)}$ formed by all $\beta\in
\mathbb{F}_q^{(2\nu)}$ such that $\alpha\, K^t\!\beta=0$ for all
$\alpha\in V$ by $V^{\bot}$. Then $[\alpha]\sim[\beta]$ if and
only if $\beta\not\in [\alpha]^{\bot}$.

 Denote the vertex set of the graph
$Sp(2\nu,q)$ by $V(Sp(2\nu,q))$. We first show that $Sp(2\nu,q)$
is strongly regular.

\begin{Theorem}
\label{parameter} $Sp(2\nu,q)$ is a strongly regular graph with
parameters
$$
\left(\frac{q^{2\nu}-1}{q-1},q^{2\nu-1},q^{2\nu-2}(q-1),q^{2\nu-2}(q-1)\right)
$$
and eigenvalues $q^{2\nu-1},q^{\nu-1}$ and $-q^{\nu-1}$.
\end{Theorem}

\begin{pf} As $|\,\mathbb{F}_q^{(2\nu)}|=q^{2\nu}$, it follows that
$|V(Sp(2\nu,q))|=\frac{q^{2\nu}-1}{q-1}$. For any $[\alpha]\in
V(Sp(2\nu,q))$, since $\dim([\alpha]^{\bot})=2\nu-1$, we see that
the degree of $[\alpha]$ which is just the number of one
dimensional subspaces $[\beta]$ such that $\beta\not\in
[\alpha]^{\bot}$, is $\frac{q^{2\nu}-q^{2\nu-1}}{q-1}=q^{2\nu-1}$.

Let $[\alpha],[\beta]$ be any two different vertices of
$Sp(2\nu,q)$ which are adjacent with each other or not. Then
$\dim([\alpha,\beta]^{\bot})=2\nu-2$. Note that a vertex
$[\gamma]$ is adjacent with both $[\alpha]$ and $[\beta]$  is
equivalent to that $\gamma\not\in
[\alpha]^{\bot}\cup[\beta]^{\bot}$. But
$$
|[\alpha]^{\bot}\cup[\beta]^{\bot}|=|[\alpha]^{\bot}|+|[\beta]^{\bot}|-|[\alpha,\beta]^{\bot}|.
$$
Hence the number of vertices which are adjacent with both
$[\alpha]$ and $[\beta]$ is
$\frac{q^{2\nu}-2q^{2\nu-1}+q^{2\nu-2}}{q-1}$ $=q^{2\nu-2}(q-1)$.
Therefore $Sp(2\nu,q)$ is a strongly regular graph with parameter
$$
\left(\frac{q^{2\nu}-1}{q-1},q^{2\nu-1},q^{2\nu-2}(q-1),q^{2\nu-2}(q-1)\right).
$$

By the same arguments as in [3, Section 10.2], we get that the
eigenvalues of $Sp(2\nu,q)$ are $q^{2\nu-1},q^{\nu-1}$ and $
-q^{\nu-1}$.
\end{pf}

Let $n\geq 2$. We say that a graph $X$ is {\em $n$-partite} if
there are subsets $X_1,\ldots, X_n$ of the vertex set $V(X)$ of
$X$ such that $V(X)=X_1\cup\cdots\cup X_n$, where $X_i\cap
X_j=\emptyset$ for all $i\not= j$, and that there is no edge of
$X$ joining two vertices of the same subset. We are going to show
that $Sp(2\nu,q)$ is $(q^{\nu}+1)$-partite. We need some results
about subspaces of $\mathbb{F}_q^{(2\nu)}$. A subspace $V$ of
$\mathbb{F}_q^{(2\nu)}$ is called {\em totally isotropic} if
$V\subseteq V^{\bot}$. Then totally isotropic subspaces of
$\mathbb{F}_q^{(2\nu)}$ are of dimension $\leq\nu$ and there exist
totally isotropic subspaces of dimension $\nu$ which are called
{\em maximal totally isotropic subspaces}, cf. \cite[Corollary
3.8]{W}.



\medskip
The following lemma is due to Dye[1].

\begin{Lemma}
\label{union} There exist maximal totally isotropic subspaces
$V_i$, $i=1,\ldots,q^{\nu}+1$, of $\mathbb{F}_q^{(2\nu)}$ such
that
$$
\mathbb{F}_q^{(2\nu)}=V_1\cup\cdots\cup V_{q^{\nu}+1},
$$
where $V_i\cap V_j=\{0\}$ for all $i\not= j$.
\end{Lemma}

\begin{Proposition}
\label{partite}
 $Sp(2\nu,q)$ is $(q^{\nu}+1)$-partite. That is,
  there exist subsets $X_1, \ldots,$ $ X_{q^\nu+1}$ of $V(Sp(2\nu,q))$ such that
$$
V(Sp(2\nu,q))=X_1\cup\cdots\cup X_{q^{\nu}+1},
$$
where $X_i\cap X_j=\emptyset$ for all $i\not=j$, and there is no
edge of $Sp(2\nu,q)$ joining two vertices of the same subset.
Moreover, the subsets $X_1, \ldots, X_{q^\nu+1}$ can be so chosen
that for any two disinct indices $i$ and $j$, every $\alpha\in
X_i$ is adjacent with exactly $q^{\nu-1}$ vertices in $X_j$.
\end{Proposition}

\begin{pf} Let $\mathbb{F}_q^{(2\nu)}=V_1\cup\cdots\cup
V_{q^{\nu}+1}$ as in \ref{union}. Set
$X_i=\{[\alpha]:\alpha\not=0\in V_i\}$, $i=1,\ldots,q^{\nu}+1$.
Then
$$
V(Sp(2\nu,q))=X_1\cup\cdots\cup X_{q^{\nu}+1}, \, X_i\cap
X_j=\emptyset, \mbox{ for all } \,i\not= j.
$$
As $V_i$ is totally isotropic, we see that there is no edge
joining any two vertices in $X_i$. Thus $Sp(2\nu,q)$ is
$(q^{\nu}+1)$-partite. For any $i\not=j$, let $[\alpha]\in X_i$.
Since $V_j$ is maximal totally isotropic of dimension $\nu$, it
follows that $\alpha\not\in V_j=V_j^{\bot}$ and
$\dim([\alpha]^{\bot}\cap V_j)=\dim([\alpha,V_j]^{\bot})=\nu-1$.
Note that, for any $[\beta]\in X_j$, $[\beta]$ is adjacent with
$[\alpha]$ if and only if $\beta\in
V_j\setminus([\alpha]^{\bot}\cap V_j)$. Hence the number of
vertices in $X_j$ which is adjacent with $[\alpha]$ is
$\frac{q^{\nu}-1}{q-1}-\frac{q^{\nu-1}-1}{q-1}=q^{\nu-1}$.
\end{pf}

Now we can compute the chromatic number of $Sp(2\nu,q)$.

\begin{Theorem}
\label{chromatic} $\chi(Sp(2\nu,q))=q^{\nu}+1$.
\end{Theorem}

\begin{pf}
By \ref{partite}, we see that $\chi(Sp(2\nu,q))\leq q^{\nu}+1$.
Note that $\chi(Sp(2\nu,q))$ is the minimal $n$ such that
$Sp(2\nu,q)$ is $n$-partite. Suppose that $Sp(2\nu,q)$ is
$n$-partite. Then there exist subsets $Y_1, \ldots, Y_n$ of
$V(Sp(2\nu,q))$ such that
$$
V(Sp(2\nu,q))=Y_1\cup\cdots\cup Y_n,\, Y_i\cap
Y_j=\emptyset,\,\mbox{ for all } \, i\not= j,
$$
and there is no edge joining any two vertices in the same $Y_i$
for $i=1,\ldots, n$. We want to show that $n\geq q^{\nu}+1$.
Suppose that $n<q^{\nu}+1$. From the above equality, we have
$\sum_{i=1}^n|Y_i|=\frac{q^{2\nu}-1}{q-1}=(\frac{q^{\nu}-1}{q-1})(q^{\nu}+1)$.
Then there exists some $i$ such that
$|Y_i|>\frac{q^{\nu}-1}{q-1}$. Let $W_i$ be the subspace of
$\mathbb{F}_q^{(2\nu)}$ generated by all $\alpha$ such that
$[\alpha]\in Y_i$. Then $W_i$ is a totally isotropic subspace,
hence $\dim W_i\leq\nu$. This turns out
$|Y_i|\leq\frac{q^{\nu}-1}{q-1}$, a contradiction. Hence
$\chi(Sp(2\nu,q))=q^{\nu}+1$.
\end{pf}

\section{Automorphisms of Symplectic Graphs}
 We recall that a $2\nu\times2\nu$ matrix $T$ is called a
 {\em symplectic matrix} (or {\em generalized symplectic matrix}) of order $2\nu$
 over $\mathbb{F}_q$ if $TK^tT=K$ (or $TK^tT=kK$ for
 some $k\in \mathbb{F}_q^*$, respectively). The set
 of symplectic matrices (or generalized symplectic matrices) of order $2\nu$
 over $\mathbb{F}_q$ forms a group with
 respect to the matrix multiplication, which is called the
 {\em symplectic group} (or {\em generalized symplectic group}, respectively) of degree $2\nu$
 over $\mathbb{F}_q$ and denoted by
 $Sp_{2\nu}(\mathbb{F}_q)$ (or $GSp_{2\nu}(\mathbb{F}_q)$). The center of $Sp_{2\nu}(\mathbb{F}_q)$
 consists of the identity
 matrix $E$
 and $-E$, and the factor group  $Sp_{2\nu}(\mathbb{F}_q)/\{E,-E\}$ is
 called the {\em projective symplectic group} of degree $2\nu$ over $\mathbb{F}_q$ and denoted by
 $PSp_{2\nu}(\mathbb{F}_q)$. The center of $GSp_{2\nu}(\mathbb{F}_q)$ consists of all
 $kE$, where $k\in \mathbb{F}_q^*$, and the factor group of  $GSp_{2\nu}(\mathbb{F}_q)$ with respect to its center
 is called the {\em projective generalized symplectic group} of degree
$2\nu$ over $\mathbb{F}_q$ and denoted by
 $PGSp_{2\nu}(\mathbb{F}_q)$. Clearly, $PGSp_{2\nu}(\mathbb{F}_q)\cong
 PSp_{2\nu}(\mathbb{F}_q)$, and when $q=2$,
 $GSp_{2\nu}(\mathbb{F}_2)=Sp_{2\nu}(\mathbb{F}_2)$.

\begin{Proposition}
\label{sympmatrix} Let $T$ be a $2\nu\times 2\nu$ nonsingular
matrix over $\mathbb{F}_q$ and
\begin{eqnarray*}
\sigma_T:\,\,V(Sp(2\nu,q))&\rightarrow&V(Sp(2\nu,q))\\
\,[\alpha]&\mapsto&[\alpha T].
\end{eqnarray*}
Then
\begin{itemize}
\item[(1)] $T\in GSp_{2\nu}(\mathbb{F}_q)$ if and only if
$\sigma_T\in\mbox{\em Aut}(Sp(2\nu,q))$. In particular, when
$q=2$, $T\in Sp_{2\nu}(\mathbb{F}_2)$ if and only if
$\sigma_T\in\mbox{\em Aut}(Sp(2\nu,2))$
\item[(2)] For any $T_1,T_2\in GSp_{2\nu}(\mathbb{F}_q)$,
$\sigma_{T_1}=\sigma_{T_2}$ if and only if $T_1=k T_2$ for some
$k\in \mathbb{F}_q$;
\end{itemize}
\end{Proposition}

\begin{pf}
It is clear that $\sigma_T$ is an one-one correspondence from
$V(Sp(2\nu,q))$ to itself.

(1) First assume $T\in GSp_{2\nu}(\mathbb{F}_q)$. Then $TK^tT=kK$
for
 some $k\in \mathbb{F}_q^*$. For any
$[\alpha],[\beta]\in V(Sp(2\nu,q))$, since $\alpha
K^t\!\beta=k^{-1}(\alpha T)K^t\!(\beta T)$, $[\alpha]\sim[\beta]$
if and only if $\sigma_T([\alpha])\sim\sigma_T([\beta])$, hence
$\sigma_T\in\mbox{Aut}(Sp(2\nu,q))$.

 Conversely, assume $\sigma_T\in\mbox{Aut}(Sp(2\nu,q))$. Then, for any
$\alpha,\beta\not=0\in \mathbb{F}_q^{(2\nu)}$, $\alpha
K^t\!\beta=0$ if and only if $\alpha(TK^tT)^t\!\beta=0$. Hence,
for any $\alpha\not=0\in \mathbb{F}_q^{(2\nu)}$, the two systems
of linear equations $(\alpha K)^t\!X=0$, $(\alpha TK^tT)^t\!X=0$
have the same solutions. But $\mbox{rank}(\alpha
K)=\mbox{rank}(\alpha TK^tT)=1$, we see that $\alpha K=k(\alpha
TK^tT)$ for some $k\in \mathbb{F}_q^*$, which depends on $\alpha$.
Take $\alpha=(1,0,\ldots,0)$,$(0,1,\ldots,0)$,
$\ldots,(0,0,\ldots,1)$, we get that
$K=\mbox{diag}(k_1,k_2,\ldots,k_{2\nu})TK^tT$, for some
$k_1,k_2,\ldots,k_{2\nu}\in \mathbb{F}_q^*$. Take
$\alpha=(1,1,\ldots,1)$, we see that $k_1=k_2=\ldots=k_{2\nu}$,
hence $K=k_1TK^tT$.

(2) It is clear that $\sigma_{T_1}=\sigma_{T_2}$ if $T_1=kT_2$ for
some $k\in \mathbb{F}_q^*$. Conversely, suppose that
$\sigma_{T_1}=\sigma_{T_2}$. Then, for any $\alpha\not=0\in
\mathbb{F}_q^{(2\nu)}$, $\alpha T_1=k\alpha T_2$ for some $k\in
\mathbb{F}_q^*$. Take $\alpha=(1,0,\ldots,0)$,$(0,1,\ldots,0)$,
and so on as above, we see that $T_1=kT_2$ for some $k\in
\mathbb{F}_q^*$.
\end{pf}


By \ref{sympmatrix}, every generalized symplectic matrix in
$GSp_{2\nu}(\mathbb{F}_q)$ induces an automorphism of $Sp(2\nu,q)$
and two generalized symplectic matrices $T_1$ and $T_2$ induce the
same automorphism of $Sp(2\nu,q)$ if and only if $T_1=kT_2$ for
some $k\in \mathbb{F}_q$. Thus $PSp_{2\nu}(\mathbb{F}_q)$ can be
regarded as a subgroup of $\mbox{Aut}(Sp(2\nu,q))$.

\begin{Proposition}
\label{transitive} $Sp(2\nu,q)$ is vertex transitive and edge
transitive.
\end{Proposition}

\begin{pf}
For any $[\alpha],[\beta]\in V(Sp(2\nu,q))$, there exists $T\in
Sp_{2\nu}(\mathbb{F}_q)$ such that $\alpha T=\beta$ by \cite[Lemma
3.11]{W}. Then $\sigma_T\in\mbox{Aut}(Sp(2\nu,q))$ such that
$\sigma_T([\alpha])=[\beta]$. Hence $Sp(2\nu,q)$ is vertex
transitive.

Let $[\alpha_1],[\alpha_2],[\beta_1],[\beta_2]\in V(Sp(2\nu,q))$
such that $[\alpha_1]\sim [\alpha_2]$ and
$[\beta_1]\sim[\beta_2]$. We may assume that
$\alpha_1K^t\!\alpha_2=\beta_1K^t\!\beta_2$. Then, by \cite[Lemma
3.11]{W} again, there exists $T\in Sp_{2\nu}(\mathbb{F}_q)$ such
that $\alpha_1T=\beta_1$ and $\alpha_2T=\beta_2$. Then
$\sigma_T\in\mbox{Aut}(Sp(2\nu,q))$ such that
$\sigma_T([\alpha_1])=[\beta_1]$ and
$\sigma_T([\alpha_2])=[\beta_2]$. Hence $Sp(2\nu,q)$ is edge
transitive.
\end{pf}

 When $q=2$, we have the following

\begin{Proposition}
$\mbox{\em Aut}(Sp(2\nu,2))\cong Sp_{2\nu}(\mathbb{F}_2)$.
\end{Proposition}

\begin{pf}
Let
\begin{eqnarray*}
\sigma:\,\, Sp_{2\nu}(\mathbb{F}_2) &\rightarrow&
\mbox{Aut}(Sp(2\nu,2))\\
T&\mapsto&\sigma_T.
\end{eqnarray*}
Then, by \ref{sympmatrix}, $\sigma$ is an injection. Clearly,
$\sigma$ preserves the operation. It remains to show that, for any
$\tau\in\mbox{Aut}(Sp(2\nu,2))$, there exists a $T\in
Sp_{2\nu}(\mathbb{F}_2)$ such that $\tau=\sigma_T$.

Note that, for any $\alpha\not=0\in \mathbb{F}_2^{(2\nu)}$, we
have that $[\alpha]=\{0,\alpha\}$. We will denote the uniquely
defined element $\tau([\alpha])\setminus\{0\}$ by $\tau(\alpha)$
and set $\tau(0)=0$. Then from $\tau\in\mbox{Aut}(Sp(2\nu,2))$ we
see that $\alpha K^t\!\beta=\tau(\alpha)K^t\!(\tau(\beta))$ for
any $\alpha,\beta\in \mathbb{F}_2^{(2\nu)}$ (not necessarily
non-zero). Fix any $\alpha\in \mathbb{F}_2^{(2\nu)}$. Let
$\beta_1,\beta_2\in \mathbb{F}_2^{(2\nu)}$. Then
\begin{eqnarray*}
\alpha K^t\!\beta_1&=&\tau(\alpha)K^t\!(\tau(\beta_1)),\\
\alpha K^t\!\beta_2&=&\tau(\alpha)K^t\!(\tau(\beta_2)).
\end{eqnarray*}
Thus $$\alpha
K^t\!(\beta_1+\beta_2)=\tau(\alpha)K^t\!(\tau(\beta_1)+\tau(\beta_2)).$$
But $$\alpha
K^t\!(\beta_1+\beta_2)=\tau(\alpha)K^t\!(\tau(\beta_1+\beta_2)),$$
hence $$\tau(\alpha)
K^t\!(\tau(\beta_1+\beta_2)+\tau(\beta_1)+\tau(\beta_2))=0.$$ This
is true for any $\alpha\in \mathbb{F}_2^{(2\nu)}$, it follows that
$\tau(\beta_1+\beta_2)+\tau(\beta_1)+\tau(\beta_2)=0$, i.e.,
$\tau(\beta_1+\beta_2)=\tau(\beta_1)+\tau(\beta_2)$. Set
$$
T=\left( \begin{array}{c} \tau(1,0,\ldots,0)\\
\tau(0,1,\ldots,0)\\
\vdots\\
\tau(0,0,\ldots,1)
\end{array}
\right).
$$
Then $\tau(\alpha)=\alpha T$ for any $\alpha\in
\mathbb{F}_2^{(2\nu)}$. Thus $T$ is nonsingular. By
\ref{sympmatrix} $T\in Sp_{2\nu}(\mathbb{F}_2)$ and
$\tau=\sigma_T$ as required.
\end{pf}

From now on, we assume that $q>2$. In $\mathbb{F}_q^{(2\nu)}$, let
us set
\begin{eqnarray*}
e_1&=&(1,0,0,0,\ldots,0,0),\\
f_1&=&(0,1,0,0,\ldots,0,0),\\
e_2&=&(0,0,1,0,\ldots,0,0),\\
f_2&=&(0,0,0,1,\ldots,0,0),\\
&&\ldots\ldots\\
e_{\nu}&=&(0,0,0,0,\ldots,1,0),\\
f_{\nu}&=&(0,0,0,0,\ldots,0,1).
\end{eqnarray*}
Then $e_i,f_i,\,i=1,\ldots,\nu$, form a basis of
$\mathbb{F}_q^{(2\nu)}$ and $e_iK^t\!f_i=1$, $e_iK^t\!e_j=0$,
$f_iK^t\!f_j=0$, $i,j=1,\ldots,\nu$, and $e_iK^t\!f_j=0$,
$i\not=j$, $i,j=1,\ldots,\nu$.

\medskip
In order to describe $\mbox{Aut}(Sp(2\nu,q))$ for any prime power
$q$, we need some definition from group theory. Let $\varphi$ be
the natural action of $\mbox{Aut}(\mathbb{F}_q)$ on the group
$\mathbb{F}_q^*\times\cdots\times \mathbb{F}_q^*$ ($\nu$ in
number) defined by
$$
\varphi(\pi)((k_1,\ldots,k_\nu))=(\pi(k_1),\ldots,\pi(k_\nu)),
\mbox{ for all }\pi\in \mbox{Aut}(\mathbb{F}_q) \mbox{ and }
k_1,\ldots,k_\nu \in \mathbb{F}_q^*,
$$ then the
semi-direct product of $\mathbb{F}_q^*\times\cdots\times
\mathbb{F}_q^*$ by $\mbox{Aut}(\mathbb{F}_q)$ corresponding to
$\varphi$, denoted by $(\mathbb{F}_q^*\times\cdots\times
 \mathbb{F}_q^*)\rtimes_{\varphi}\mbox{Aut}(\mathbb{F}_q)$, is the group consisting of all elements of
the form $(k_1,\ldots,k_\nu, \pi)$, where $k_1,\ldots,k_\nu \in
\mathbb{F}_q^*$ and $\pi\in \mbox{Aut}(\mathbb{F}_q)$, with
multiplication defined by
$$
(k_1,\ldots,k_\nu, \pi)(k_1^{'},\ldots,k_\nu^{'},\pi^{'})=(k_1\pi
(k_1^{'}),\ldots, k_\nu\pi(k_\nu^{'}),\pi\pi^{'}).$$ Then the main
result about $\mbox{Aut}(Sp(2\nu,q))$ is as follows.

\begin{Theorem}
Regard $PSp_{2\nu}(\mathbb{F}_q)$ as a subgroup of $\mbox{\em
Aut}(Sp(2\nu,q))$ and let  $E$ be the subgroup of $\mbox{\em
Aut}(Sp(2\nu,q))$ defined as follows
$$
E=\{\sigma\in \mbox{\em Aut}(Sp(2\nu,q)):
\sigma([e_i])=[e_i],\sigma([f_i])=[f_i], i=1,\ldots,\nu\}.
$$
Then
\begin{itemize}
\item[(1)] $\mbox{\em Aut}(Sp(2\nu,q))=PSp_{2\nu}(\mathbb{F}_q)\cdot E$;
\item[(2)] If $\nu=1$, then $E$ is isomorphic to the symmetric group on $q-1$ elements;
\item[(3)] If $\nu > 1$, then
$$
E\cong\underbrace{(\mathbb{F}_q^*\times\cdots\times
\mathbb{F}_q^*)}_{\nu}\rtimes_{\varphi}\mbox{\em
Aut}(\mathbb{F}_q).
$$
\end{itemize}
\end{Theorem}

\begin{pf} (1) Let $\tau\in\mbox{Aut}(Sp(2\nu,q))$. Suppose that $\tau([e_i])=[e_i']$,
$\tau([f_i])=[f_i']$, $i=1,\ldots,\nu$. Then
$e_i'K^t\!f_i'\not=0$, $e_i'K^t\!e_j'=0$, $f_i'K^t\!f_j'=0$,
$i,j=1,\ldots,\nu$ and $e_i'K^t\!f_j'=0$,  $i\not=j$,
$i,j=1,\ldots,\nu$. We may choose $e_i', f_i',\,i=1,\ldots,\nu$,
such that $e_i'K^t\!f_i'=1$,$i=1,\ldots,\nu$. Let
$$
A=\left(\begin{array}{c} e_1\\
f_1\\e_2\\f_2\\
\vdots\\
e_{\nu}\\
f_{\nu}
\end{array}\right),\,\,A'=\left(\begin{array}{c} e_1'\\
f_1'\\e_2'\\f_2'\\
\vdots\\
e_{\nu}'\\
f_{\nu}'
\end{array}\right).
$$
Then $AK^t\!A=K=A'K^t\!A'$. Thus, by \cite[Lemma 3.11]{W}, there
exists $T\in Sp_{2\nu}(\mathbb{F}_q)$ such that $A=A'T$, i.e.,
$e_i'T=e_i$, $f_i'T=f_i$, $i=1,\ldots,\nu$. Set
$\tau_1=\sigma_{T}\tau$. Then $\tau_1([e_i])=[e_i]$,
$\tau_1([f_i])=[f_i]$, $i=1,\ldots,\nu$, hence $\tau_1\in E$. Thus
$\tau\in PSp_{2\nu}(\mathbb{F}_q)\cdot E$. It follows that
$\mbox{Aut}(Sp(2\nu,q))=PSp_{2\nu}(\mathbb{F}_q)\cdot E$.

(2) When $\nu=1$, it is clear that $E$ is isomorphic to the
symmetric group on the $q-1$ vertices of $Sp(2,q)$ since $Sp(2,q)$
is a complete graph.

(3) Suppose that $\nu>1$. Our idea to prove the third part of the
theorem is to identify some elements of $E$ which form a subgroup
of $E$ isomorphic to the semi-direct product in the theorem and,
then, to show that every element of $E$ has the form of the
elements identified before.

Firstly, let us write out some elements of $E$. Let
$k_1,\ldots,k_{\nu}\in \mathbb{F}^*_q$ and
$\pi\in\mbox{Aut}(\mathbb{F}_q)$. Let
$\sigma_{(k_1,\ldots,k_{\nu},\pi)}$ be the map which takes any
vertex $[a_1,a_2,a_3,a_4,\ldots,a_{2\nu-1},a_{2\nu}]$ of
$Sp(2\nu,q)$ to the vertex
$$
[\pi(a_1),k_1\pi(a_2),k_2\pi(a_3),k_1k_2^{-1}\pi(a_4),\ldots,k_{\nu}
\pi(a_{2\nu-1}),k_1k_{\nu}^{-1}\pi(a_{2\nu})].
$$
Then it is clear that $\sigma_{(k_1,\ldots,k_{\nu},\pi)}$ is
well-defined. Furthermore, it is easy to see that
$\sigma_{(k_1,\ldots,k_{\nu},\pi)}$ is injective, but the vertex
set of $Sp(2\nu,q)$ is finite, $\sigma_{(k_1,\ldots,k_{\nu},\pi)}$
is a bijection from $V(Sp(2\nu,q))$ to itself. Let
$\alpha=[a_1,a_2,a_3,a_4,\ldots,a_{2\nu-1},a_{2\nu}]$,
$\beta=[a'_1,a'_2,a'_3,a'_4,\ldots,a'_{2\nu-1},a'_{2\nu}]$ be two
vertices of $Sp(2\nu,q)$. If $\alpha\not\sim\beta$, then, by
definition,
$$
(a_1a'_2-a_2a'_1)+(a_3a'_4-a_4a'_3)+\ldots+(a_{2\nu-1}a'_{2\nu}-a_{2\nu}a'_{2\nu-1})=0,
$$
which implies that \begin{eqnarray*} &&
(\pi(a_1)k_1\pi(a'_2)-\pi(a_2)k_1\pi(a'_1))+(k_2\pi(a_3)k_1k_2^{-1}\pi(a'_4)
-k_1k_2^{-1}\pi(a_4)k_2\pi(a'_3))\\
&&+\ldots+(k_{\nu}\pi(a_{2\nu-1})k_1k_{\nu}^{-1}\pi(a'_{2\nu})
-k_1k_{\nu}^{-1}\pi(a_{2\nu})k_{\nu}\pi(a'_{2\nu-1}))=0,
\end{eqnarray*}
i.e.,
$\sigma_{(k_1,\ldots,k_{\nu},\pi)}(\alpha)\not\sim\sigma_{(k_1,\ldots,k_{\nu},\pi)}(\beta)$.
Since the edges set of $Sp(2\nu,q)$ is finite,
$\alpha\not\sim\beta$ if and only if
$\sigma_{(k_1,\ldots,k_{\nu},\pi)}(\alpha)\not\sim\sigma_{(k_1,\ldots,k_{\nu},\pi)}(\beta)$.
Hence
$\sigma_{(k_1,\ldots,k_{\nu},\pi)}\in\mbox{Aut}(Sp(2\nu,q))$. Note
that $\sigma_{(k_1,\ldots,k_{\nu},\pi)}([e_i])=[e_i]$, $
\sigma_{(k_1,\ldots,k_{\nu},\pi)}([f_i])=[f_i]$, $i=1,\ldots,\nu$,
hence, $\sigma_{(k_1,\ldots,k_{\nu},\pi)}\in E$.

If we define a map $h$ as
$(k_1,\ldots,k_{\nu},\pi)\mapsto\sigma_{(k_1,\ldots,k_{\nu},\pi)}$,
then it is easy to verify that $h$ is a group homomorphism from
$(\mathbb{F}_q^*\times\cdots\times
\mathbb{F}_q^*)\times_{\varphi}\mbox{Aut}(\mathbb{F}_q)$ to $E$.
It is also easy to see that if $(k_1,\ldots,
k_{\nu},\pi)\not=(k'_1,\ldots, k'_{\nu},\pi')$ then
$\sigma_{(k_1,\ldots,k_{\nu},\pi)}\not=
\sigma_{(k'_1,\ldots,k'_{\nu},\pi')}$. In order to prove the third
part of the theorem, we will show that $h$ is a group isomorphism.
It remains to show that every element of $E$ is of the form
$\sigma_{(k_1,\ldots,k_{\nu},\pi)}$.

Suppose that $\sigma\in E$. Note that if
$\sigma([a_1,a_2,\ldots,a_{2\nu}])=[b_1,b_2,\ldots,b_{2\nu}]$,
then $a_{2i-1}\not=0$ if and only if
$[a_1,a_2,\ldots,a_{2\nu}]\sim [f_i]$ and $a_{2i}\not=0$ if and
only if $[a_1,a_2,\ldots,a_{2\nu}]$ $\sim [e_i]$, and similar
results are also true for $b_i$. But $\sigma([e_i])=[e_i]$ and
$\sigma([f_i])=[f_i]$, it follows that $a_i=0$ if and only if
$b_i=0$. For any vertex $[a_1,a_2,\ldots,a_{2\nu}]$, if
$a_1=\cdots=a_{i-1}=0$ and $a_i\not=0$ then
$[a_1,a_2,\ldots,a_{2\nu}]$ can be uniquely written as
$[0,\ldots,0,1,a'_{i+1},\ldots,a'_{2\nu}]$ and
$\sigma([a_1,a_2,\ldots,a_{2\nu}])$ can be uniquely written as
$[0,\ldots,0,1,b'_{i+1},\ldots,b'_{2\nu}]$. Let us show how to
determine $b'_{i+1},\ldots,b'_{2\nu}$ from $a'_{i+1},\ldots,$ $
a'_{2\nu}$. We will use frequently the fact that, for any vertices
$[\alpha],[\beta]$, if $[\alpha]\not\sim[\beta]$ then
$\sigma([\alpha])\not\sim\sigma([\beta])$.

In the following, we will denote
$[a_1,a_1',a_2,a'_2,\ldots,a_{\nu},a'_{\nu}]$ by
$\sum_{i=1}^{\nu}a_i[e_i]+\sum_{i=1}^{\nu}a'_i[f_i]$, for example,
$[a,b,0,\ldots,0]$ is denoted by $a[e_1]+b[f_1]$. Since $\sigma$
is a bijection from $V(Sp(2\nu,q))$ to itself, we have
permutations $\pi_i$, $i=2,\ldots,2\nu$, of $\mathbb{F}_q$ with
$\pi(0)=0$ such that
\begin{eqnarray*}
\sigma([e_1]+a_{2i-1}[e_i])&=&[e_1]+\pi_{2i-1}(a_{2i-1})[e_i],\\
\sigma([e_1]+a_{2i}[f_i])&=&[e_1]+\pi_{2i}(a_{2i})[f_i].
\end{eqnarray*}

We firstly consider the cases $\sigma([0,1,a_3,\ldots,a_{2\nu}])$
and $\sigma([1,a_2,a_3,\ldots,a_{2\nu}])$. Let
$\sigma([0,1,a_3,\ldots,a_{2\nu}])=[0,1,a_3',\ldots,a'_{2\nu}]$
and $j\geq 1$. If $a_{2j+1}\not=0$, then, from\\
$[0,1,a_3,\ldots,a_{2\nu}]\not\sim [e_1]+a_{2j+1}^{-1}[f_{j+1}]$
we have $[0,1,a'_3,\ldots,a'_{2\nu}]\not\sim
[e_1]+\pi_{2j+2}(a_{2j+1}^{-1})[f_{j+1}]$, hence,
$a'_{2j+1}=\pi_{2j+2}(a_{2j+1}^{-1})^{-1}$. If $a_{2j+2}\not=0$,
then from $[0,1,a_3,\ldots,a_{2\nu}]\not\sim
[e_1]-a_{2j+2}^{-1}[e_{j+1}]$ we have
$[0,1,a'_3,\ldots,a'_{2\nu}]\not\sim
[e_1]+\pi_{2j+1}(-a_{2j+2}^{-1})[e_{j+1}]$, hence,
$a'_{2j+2}=-\pi_{2j+1}(-a_{2j+2}^{-1})^{-1}$. Thus
\begin{eqnarray}
\sigma([0,1,a_3,\ldots,a_{2\nu}])=[0,1,a_3',\ldots,a'_{2\nu}],
\end{eqnarray}
where $a'_{2j+1}=\pi_{2j+2}(a_{2j+1}^{-1})^{-1}$ if
$a_{2j+1}\not=0$ and $a'_{2j+2}=-\pi_{2j+1}(-a_{2j+2}^{-1})^{-1}$
if $a_{2j+2}\not=0$.

For the case $\sigma([1,a_2,a_3,\ldots,a_{2\nu}])$. Let
$\sigma([1,a_2,a_3,\ldots,a_{2\nu}])=[1,a''_2,a''_3,\ldots,a''_{2\nu}]$.
From $[1,a_2,a_3,\ldots,a_{2\nu}]\not\sim [e_1]+a_2[f_1]$ we get
$[1,a''_2,a''_3,\ldots,a''_{2\nu}]\not\sim [e_1]+\pi_2(a_2)[f_1]$,
hence, $a''_2=\pi_2(a_2)$. Let $j\geq 1$. If $a_{2j+1}\not=0$,
then, from $[1,a_2,a_3,\ldots,a_{2\nu}]\not\sim
[f_1]-a_{2j+1}^{-1}[f_{j+1}]$ and
$\sigma([f_1]-a_{2j+1}^{-1}[f_{j+1}])=[f_1]-\pi_{2j+1}(a_{2j+1})^{-1}[f_{j+1}]$
as shown above, we have $[1,a''_2,a''_3,\ldots,a''_{2\nu}]\not\sim
[f_1]-\pi_{2j+1}(a_{2j+1})^{-1}[f_{j+1}]$, hence,
$a''_{2j+1}=\pi_{2j+1}(a_{2j+1})$. Similarly, if $a_{2j+2}\not=0$,
then from $[1,a_2,a_3,\ldots,a_{2\nu}]\not\sim
[f_1]+a_{2j+2}^{-1}[e_{j+1}]$ we have
$[1,a''_2,a''_3,\ldots,a''_{2\nu}]\not\sim
[f_1]+\pi_{2j+2}(a_{2j+2})^{-1}[e_{j+1}]$, hence,
$a''_{2j+2}=\pi_{2j+2}(a_{2j+2})$. Thus, for any
$a_2,a_3,\ldots,a_{2\nu}\in \mathbb{F}_q$,
\begin{eqnarray}
\sigma([1,a_2,a_3,\ldots,a_{2\nu}])
=[1,\pi_2(a_2),\pi_3(a_3),\ldots,\pi_{2\nu}(a_{2\nu})].
\end{eqnarray}

Then, let $i\geq 2$, we discuss the general cases
$\sigma([0,\ldots,0,1,a_{2i+1},\ldots,a_{2\nu}])$ and
$\sigma([0,\ldots,0,1,a_{2i},\ldots,a_{2\nu}])$. The above results
of case $i=1$ will be used. Let
$\sigma([0,\ldots,0,1,a_{2i+1},\ldots,a_{2\nu}])=[0,\ldots,0,1,a'_{2i+1},\ldots,a'_{2\nu}]$
and $j\geq i$. If $a_{2j+1}\not=0$, then, from
$$
[0,\ldots,0,1,a_{2i+1},\ldots,a_{2\nu}]\not\sim
[e_1]+[e_{i}]+a_{2j+1}^{-1}[f_{j+1}]
$$
and
$\sigma([e_1]+[e_{i}]+a_{2j+1}^{-1}[f_{j+1}])=[e_1]+\pi_{2i-1}(1)[e_{i}]+\pi_{2j+2}(a_{2j+1}^{-1})[f_{j+1}]$
as shown above, we have
$$
[0,\ldots,0,1,a'_{2i+1},\ldots,a'_{2\nu}]\not\sim
[e_1]+\pi_{2i-1}(1)[e_{i}]+\pi_{2j+2}(a_{2j+1}^{-1})[f_{j+1}],
$$
hence, $a'_{2j+1}=\pi_{2i-1}(1)\pi_{2j+2}(a_{2j+1}^{-1})^{-1}$.
Similarly, if $a_{2j+2}\not=0$, then from
$$
[0,\ldots,0,1,a_{2i+1},\ldots,a_{2\nu}]\not\sim
[e_1]+[e_{i}]-a_{2j+2}^{-1}[e_{j+1}]
$$
we have
$$
[0,\ldots,0,1,a'_{2i+1},\ldots,a'_{2\nu}]\not\sim
[e_1]+\pi_{2i-1}(1)[e_{i}]+\pi_{2j+1}(-a_{2j+2}^{-1})[e_{j+1}],
$$
hence, $a'_{2j+2}=-\pi_{2i-1}(1)\pi_{2j+1}(-a_{2j+2}^{-1})^{-1}$.
Thus,
\begin{eqnarray}
\sigma([0,\ldots,0,1,a_{2i+1},\ldots,a_{2\nu}])=[0,\ldots,0,1,a'_{2i+1},\ldots,a'_{2\nu}],
\end{eqnarray}
where $a'_{2j+1}=\pi_{2i-1}(1)\pi_{2j+2}(a_{2j+1}^{-1})^{-1}$ if
$a_{2j+1}\not=0$ and
$a'_{2j+2}=-\pi_{2i-1}(1)\pi_{2j+1}(-a_{2j+2}^{-1})^{-1}$ if
$a_{2j+2}\not=0$.

Finally, for the case
$\sigma([0,\ldots,0,1,a_{2i},\ldots,a_{2\nu}])$. Let
$\sigma([0,\ldots,0,1,a_{2i},\ldots,a_{2\nu}])=$\\ $
[0,\ldots,0,1,a''_{2i},\ldots,a''_{2\nu}]$. From
$$
[0,\ldots,0,1,a_{2i},\ldots,a_{2\nu}]\not\sim
[e_1]+[e_{i}]+a_{2i}[f_i]
$$
we get
$$
[0,\ldots,0,1,a''_{2i},\ldots,a''_{2\nu}]\not\sim
[e_1]+\pi_{2i-1}(1)[e_{i}]+\pi_{2i}(a_{2i})[f_i],
$$
hence,
$a''_{2i}=\pi_{2i-1}(1)^{-1}\pi_{2i}(a_{2i})$. Let $j\geq i$. If
$a_{2j+1}\not=0$, then from
$$
[0,\ldots,0,1,a_{2i},\ldots,a_{2\nu}]\not\sim
[f_i]-a_{2j+1}^{-1}[f_{j+1}]
$$
we have
$$
[0,\ldots,0,1,a''_{2i},\ldots,a''_{2\nu}]\not\sim
[f_i]-\pi_{2i-1}(1)^{-1}\pi_{2j+1}(a_{2j+1})^{-1}[f_{j+1}],
$$
hence,
$a''_{2j+1}=\pi_{2i-1}(1)^{-1}\pi_{2j+1}(a_{2j+1})$. If
$a_{2j+2}\not=0$, then from
$$
[0,\ldots,0,1,a_{2i},\ldots,a_{2\nu}]\not\sim
[f_i]+a_{2j+2}^{-1}[e_{j+1}]
$$
we have
$$[0,\ldots,0,1,a''_{2i},\ldots,a''_{2\nu}]\not\sim
[f_i]+\pi_{2i-1}(1)\pi_{2j+2}(a_{2j+2})^{-1}[e_{j+1}],
$$
hence, $a''_{2j+2}=\pi_{2i-1}(1)^{-1}\pi_{2j+2}(a_{2j+2})$. Thus,
for any $a_{2i},a_{2i+1},\ldots,a_{2\nu}\in \mathbb{F}_q$,
\begin{eqnarray*}
(4)\,\,&&\sigma([0,\ldots,0,1,a_{2i},a_{2i+1},\ldots,a_{2\nu}])\\
&=&[0,\ldots,0,1,\pi_{2i-1}(1)^{-1}\pi_{2i}(a_{2i}),
\pi_{2i-1}(1)^{-1}\pi_{2i+1}(a_{2i+1}),\ldots,\pi_{2i-1}(1)^{-1}\pi_{2\nu}(a_{2\nu})].
\end{eqnarray*}

\medskip
Having represented $\sigma$ by $\pi_i$, $i=2,\ldots,2\nu$, let us
discuss some properties of $\pi_i$.

\begin{Lemma}
\label{pi}
\begin{itemize}
\item[(1)] For any $i\geq 1$ and $a\in \mathbb{F}_q$,
$$
\pi_{2i+1}(1)\pi_{2i+2}(a)=\pi_{2i+2}(1)\pi_{2i+1}(a)=\pi_2(a);
$$
\item[(2)] For any $i\geq 2$ and $a,b\in \mathbb{F}_q$,
\begin{eqnarray*}
\pi_i(a+b)&=&\pi_i(a)+\pi_i(b);\\
\pi_i(-a)&=&-\pi_i(a);\\
\pi_i(ab)&=&\pi_i(a)\pi_i(b)\pi_i(1)^{-1};\\
\pi_i(a^{-1})&=&\pi_i(a)^{-1}\pi_i(1)^2 \mbox{ if $a\not=0$ }.
\end{eqnarray*}
\end{itemize}
\end{Lemma}

\begin{pf} (1) We may assume that $a\not=0$. Since
$[e_1]+a[e_{i+1}]+a[f_{i+1}]\not\sim [e_{i+1}]+[f_{i+1}]$, it
follows that $\sigma([e_1]+a[e_{i+1}]+a[f_{i+1}])\not\sim
\sigma([e_{i+1}]+[f_{i+1}])$, but
\begin{eqnarray*}
\sigma([e_1]+a[e_{i+1}]+a[f_{i+1}])&=&[e_1]+\pi_{2i+1}(a)[e_{i+1}]+\pi_{2i+2}(a)[f_{i+1}],\\
\sigma([e_{i+1}]+[f_{i+1}])&=&[e_{i+1}]+\pi_{2i+1}(1)^{-1}\pi_{2i+2}(1)[f_{i+1}],
\end{eqnarray*}
we have that
$$
\pi_{2i+1}(1)^{-1}\pi_{2i+2}(1)\pi_{2i+1}(a)-\pi_{2i+2}(a)=0,
$$
i.e.,
$$
\pi_{2i+1}(1)\pi_{2i+2}(a)=\pi_{2i+2}(1)\pi_{2i+1}(a).
$$
Similarly, since $[e_1]+a[f_1]+[e_{i+1}]\not\sim
[e_1]+a[f_{i+1}]$, we have that
$[e_1]+\pi_2(a)[f_1]+\pi_{2i+1}(1)[e_{i+1}]\not\sim
[e_1]+\pi_{2i+2}(a)[f_{i+1}]$, hence,
$\pi_{2i+1}(1)\pi_{2i+2}(a)=\pi_2(a)$.

(2) From $[e_1]+(a+b)[f_1]+[e_2]\not\sim [e_1]+a[f_1]+b[f_2]$ we
have that
$$
[e_1]+\pi_2(a+b)[f_1]+\pi_3(1)[e_2]\not\sim
[e_1]+\pi_2(a)[f_1]+\pi_4(b)[f_2].
$$
Then
$\pi_2(a)-\pi_2(a+b)+\pi_3(1)\pi_4(b)=0$, but
$\pi_3(1)\pi_4(b)=\pi_2(b)$, hence,
$\pi_2(a+b)=\pi_2(a)+\pi_2(b)$. It turns out from (1) that this
equality holds for all $i\geq 2$. Thus $\pi_i(-a)=-\pi_i(a)$ as
$\pi_i(0)=0$.

For multiplication, let $i\geq 1$, from
$[e_1]+b[e_{i+1}]+ab[f_{i+1}]\not\sim [e_{i+1}]+a[f_{i+1}]$ we get
that
$$
[e_1]+\pi_{2i+1}(b)[e_{i+1}]+\pi_{2i+2}(ab)[f_{i+1}]\not\sim
[e_{i+1}]+\pi_{2i+1}(1)^{-1}\pi_{2i+2}(a)[f_{i+1}],
$$
hence,
$\pi_{2i+1}(b)\pi_{2i+1}(1)^{-1}\pi_{2i+2}(a)-\pi_{2i+2}(ab)=0$,
but
$\pi_{2i+1}(b)\pi_{2i+1}(1)^{-1}=\pi_{2i+2}(b)\pi_{2i+2}(1)^{-1}$.
Thus
$$
\pi_{2i+2}(ab)=\pi_{2i+2}(a)\pi_{2i+2}(b)\pi_{2i+2}(1)^{-1}.
$$
It follows from
$\pi_{2i+1}(1)\pi_{2i+2}(a)=\pi_{2i+2}(1)\pi_{2i+1}(a)$ and
$\pi_{2i+1}(1)\pi_{2i+2}(1)=\pi_{2i+2}(1)\pi_{2i+1}(1)$ that the
abve equality also holds for $2i+1$. It remains to consider
$\pi_2$. We have
\begin{eqnarray*}
\pi_2(ab)&=&\pi_3(1)\pi_4(ab)\\
&=&\pi_3(1)\pi_4(1)^{-1}\pi_4(a)\pi_4(b)\\
&=&\pi_3(1)^{-1}\pi_4(1)^{-1}\pi_2(a)\pi_2(b)\\
&=&\pi_2(a)\pi_2(b)\pi_2(1)^{-1}.
\end{eqnarray*}

Finally, if $a\not=0$, then from
$\pi_i(1)=\pi_i(aa^{-1})=\pi_i(a)\pi_i(a^{-1})\pi_i(1)^{-1}$ we
obtain that $\pi_i(a^{-1})=\pi_i(a)^{-1}\pi_i(1)^2$, then the
proof of lemma is complete.
\end{pf}

We continue the proof of the theorem.  Let us denote the identity
automorphism on $\mathbb{F}_q$ by $\pi_1$. Then when $i=1$, (3)
reduces to (1) and (4) reduces to (2). Therefore (3) and (4) hold
for all $i$, where $1\leq i \leq \nu$. By the above lemma, for any
$i\geq 1$, we can rewrite (3) in the form of (4) as follows.
In (3), for any $j\geq i$, we have
\begin{eqnarray*}
a'_{2j+1}&=&\pi_{2i-1}(1)\pi_{2j+2}(a_{2j+1}^{-1})^{-1}\\
&=&\pi_{2i-1}(1)\pi_{2j+2}(a_{2j+1})\pi_{2j+2}(1)^{-2}\\
&=&\pi_{2i-1}(1)\pi_{2j+1}(1)^{-1}\pi_{2j+2}(1)^{-1}\pi_{2j+1}(a_{2j+1})\\
&=&\pi_{2i-1}(1)\pi_2(1)^{-1}\pi_{2j+1}(a_{2j+1})\\
&=&\pi_{2i}(1)^{-1}\pi_{2j+1}(a_{2j+1}),
\end{eqnarray*}
and
\begin{eqnarray*}
a'_{2j+2}&=&-\pi_{2i-1}(1)\pi_{2j+1}(-a_{2j+2}^{-1})^{-1}\\
&=&\pi_{2i-1}(1)\pi_{2j+1}(a_{2j+2}^{-1})^{-1}\\
&=&\pi_{2i-1}(1)\pi_{2j+1}(a_{2j+2})\pi_{2j+1}(1)^{-2}\\
&=&\pi_{2i-1}(1)\pi_{2j+1}(1)^{-1}\pi_{2j+2}(1)^{-1}\pi_{2j+2}(a_{2j+2})\\
&=&\pi_{2i-1}(1)\pi_2(1)^{-1}\pi_{2j+2}(a_{2j+2})\\
&=&\pi_{2i}(1)^{-1}\pi_{2j+2}(a_{2j+2}).
\end{eqnarray*}
Hence, for any $a_{2i+1},\ldots,a_{2\nu}\in \mathbb{F}_q$,
$$
(5)\,\,\sigma([0,\ldots,0,1,a_{2i+1},\ldots,a_{2\nu}])=[0,\ldots,0,1,\pi_{2i}(1)^{-1}\pi_{2i+1}(a_{2i+1}),
\ldots,\pi_{2i}(1)^{-1}\pi_{2\nu}(a_{2\nu})],
$$
which is of the same form as (4).

Now let $k_1=\pi_2(1)$, $\pi=k_1^{-1}\pi_2$, $k_2=\pi_3(1)$,
$k_3=\pi_5(1), \ldots, k_{\nu}=\pi_{2\nu-1}(1)$. Then
$\pi\in\mbox{Aut}(\mathbb{F}_q)$, $\pi_2=k_1\pi$,
$\pi_3=k_2\pi,\pi_4=k_1k_2^{-1}\pi,\ldots,\pi_{2\nu-1}=k_{\nu}\pi,\pi_{2\nu}=k_1k_{\nu}^{-1}\pi$.
Assembling (4) and (5), we obtain
\begin{eqnarray*}
&&\sigma([a_1,a_2,a_3,a_4,\ldots,a_{2\nu-1},a_{2\nu}])\\
&=&[\pi(a_1),k_1\pi(a_2),k_2\pi(a_3),k_1k_2^{-1}\pi(a_4),\ldots,k_{\nu}
\pi(a_{2\nu-1}),k_1k_{\nu}^{-1}\pi(a_{2\nu})].
\end{eqnarray*}
Hence $\sigma=h(k_1,\ldots,k_{\nu},\pi)$ as required, the proof of
the theorem is complete.
\end{pf}

\begin {Corollary} When $\nu=1$,
$$|\mbox{\em Aut}(Sp(2,q))|=q(q^2-1)\cdot(q-2)!,$$ and when $\nu\geq 2$,
$$|\mbox{\em Aut}(Sp(2\nu,q))|=q^{\nu^2}\prod^{\nu}_{i=1}(q^{2i}-1)\cdot[\mathbb{F}_q:\mathbb{F}_p].$$
\end{Corollary}

\begin{pf} Note that $PSp_{2\nu}(\mathbb{F}_q)\cap E$ consists of $\sigma$ which
is reduced from some matrix of the form
$\mbox{diag}(k_1,l_1,k_2,l_2,\ldots\ldots,k_{\nu},l_{\nu})$, with
$k_il_i=1$, $i=1,\ldots,\nu$. Thus $|PSp_{2\nu}(\mathbb{F}_q)\cap
E|=\frac{1}{2}(q-1)^{\nu}$. Hence
\begin{eqnarray*}
|\mbox{Aut}(Sp(2\nu,q))|&=&\frac{|PSp_{2\nu}(\mathbb{F}_q)|
|E|}{|PSp_{2\nu}(\mathbb{F}_q)\cap E|}\\
&=& \frac{\frac{1}{2}q^{\nu^2}\prod^{\nu}_{i=1}(q^{2i}-1)\cdot
|E|}{\frac{1}{2}(q-1)^{\nu}}.
\end{eqnarray*}
Thus, when $\nu=1$, $|\mbox{Aut}(Sp(2,q))|=q(q^2-1)\cdot(q-2)!$,
and when $\nu\geq 2$,
\begin{eqnarray*}
&&|\mbox{Aut}(Sp(2\nu,q))|\\
&=& \frac{\frac{1}{2}q^{\nu^2}\prod^{\nu}_{i=1}(q^{2i}-1)\cdot
(q-1)^{\nu}\cdot|\mbox{Aut}(\mathbb{F}_q)|}{\frac{1}{2}(q-1)^{\nu}}\\
&=&q^{\nu^2}\prod^{\nu}_{i=1}(q^{2i}-1)\cdot
|\mbox{Aut}(\mathbb{F}_q)|\\
&=&q^{\nu^2}\prod^{\nu}_{i=1}(q^{2i}-1)\cdot[\mathbb{F}_q:\mathbb{F}_p],
\end{eqnarray*}
as is well-known that
$|\mbox{Aut}(\mathbb{F}_q)|=[\mathbb{F}_q:\mathbb{F}_p]$ where
$p=\mbox{char}(\mathbb{F}_q)$.
\end{pf}

\medskip
\medskip
\noindent
Zhongming Tang\\
Department of Mathematics\\
Suzhou University\\
Suzhou 215006\\
P.\ R.\ China\\
E-mail: zmtang@suda.edu.cn

\medskip
\noindent
Zhe-xian Wan\\
Academy of Mathematics and System Sciences\\
Chinese Academy of Science\\
Beijing 100080\\
P.\ R.\ China\\
E-mail: wan@amss.ac.cn

\end{document}